\documentclass[12pt]{article}
\usepackage{graphicx}

\newtheorem{theorem}{Theorem}
\newtheorem{corollary}[theorem]{Corollary}
\newenvironment{proof}[1][Proof]{\textbf{#1.} }{\ \rule{0.5em}{0.5em}}

\def\text{\hbox} \def\newpage{\vfill\break}

\def\a{\alpha}

\def\p{\pi}

\def\s{\sigma}

\def\t{\tau}

\def\L{\Lambda}

\def\S{\Sigma}

\def\T{{\mathcal T}}
\def\K{{\mathcal K}}
\def\L{{\mathcal L}}
\def\Z{{\mathbf Z}}
\def\S{{\mathbf S}}
\def\wh{\widehat}

\def\T1{T_1({\bf p/q};{\bf r/s})}

\def\K2{{\mathcal K}_2({\bf p/q};{\bf r/s})}

\begin{document}

\title{Generalized Takahashi manifolds
\thanks{Work performed
under the auspices of G.N.S.A.G.A. of C.N.R. of Italy and
supported by the University of Bologna, funds for selected
research topics. The second named author was supported by RFBR
(grant no. 98-01-00699).}}

\author{Michele Mulazzani \and Andrei Vesnin}

\maketitle

\begin{abstract}
We introduce a family of closed 3-dimensional manifolds, which are
a generalization of certain manifolds studied by M. Takahashi. The
manifolds are represented by Dehn surgery with rational
coefficients on $\S^3$, along an $n$-periodic $2n$-component link.
A presentation of their fundamental group is obtained, and
covering properties of these ma\-ni\-folds are studied. In
particular, this family of manifolds includes the whole class of
cyclic branched coverings of two-bridge knots. As a consequence we
obtain a simple explicit surgery presentation for this important
class of manifolds.
\end{abstract}

\section{Introduction}

\begin{figure}[bht]
 \begin{center}
 \includegraphics*[totalheight=3cm]{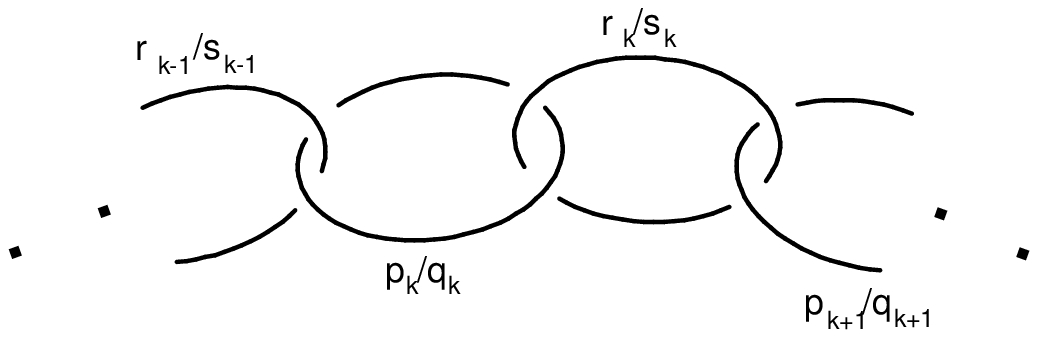}
 \end{center}
 \caption{Surgery presentation for Takahashi manifolds.}
 \label{Fig. 0}
\end{figure}

Takahashi manifolds are closed orientable 3-manifolds introduced
in \cite{Ta} by Dehn surgery with rational coefficients on ${\bf
S}^3$, along the $2n$-component link ${\cal L}_{n}$ of Figure
\ref{Fig. 0}, which is a closed chain of $2n$ unknotted components.
These manifolds have been intensively studied in
\cite{KKV}, \cite{KV}, \cite{Mu2} and \cite{RS}. In particular, a
topological characterization of all Takahashi manifolds as
two-fold coverings of ${\bf S}^3$, branched over the closure of
certain rational 3-string braids, is given in \cite{KV} and
\cite{RS}.

A Takahashi manifold is said to be {\it periodic\/} when the
surgery coefficients have the same cyclic symmetry of order $n$ of
the link ${\cal L}_{n}$, i.e. the coefficients are $p_k/q_k=p/q$
and $r_k/s_k=r/s$ alternately, for $k=1,\dots,n$. Several
important classes of 3-manifolds, such as (fractional) Fibonacci
manifolds \cite{HKM,KV} and Sieradski manifolds \cite{CHK},
represent notable examples of periodic Takahashi manifolds. More
generally, all cyclic branched coverings of two-bridge knots of
genus one are periodic Takahashi manifolds \cite{KKV}. A
characterization of periodic Takahashi manifolds as $n$-fold
cyclic coverings of the connected sum of two lens spaces, branched
over a knot, is given in \cite{Mu2}.

In this paper we generalize the family of Takahashi manifolds, as
well as periodic Takahashi manifolds, considering surgery along a
more general family of links (see Figure \ref{Fig. 1}). We obtain
a presentation for the fundamental groups (Theorem
\ref{fundamental}) and study covering properties of these
manifolds. The generalized Takahashi manifolds are described as
2-fold branched coverings of $\S^3$ (Theorem \ref{twofold}) and
the genera\-li\-zed periodic Takahashi manifolds are described as
the $n$-fold cyclic branched coverings of the connected sum of
lens spaces (Theorem \ref{cyclic}). In particular, we show that
the family of generalized periodic Takahashi manifolds contains
all cyclic coverings of two-bridge knots (Corollary \ref{all
twobridge}), thus obtaining a simple explicit surgery presentation
for this important class of manifolds (Figure \ref{Fig. 00}). This
shows that our generalization of Takahashi manifolds is, in this
sense, really natural. As a further result, we give cyclic presentations
(in the sense of \cite{Jo}) for the fundamental groups of all
cyclic branched coverings of two-bridge knots (Theorem
\ref{cyclicpresentation}).

\section{Construction of the manifolds}

\begin{figure}[bht]
 \begin{center}
 \includegraphics*[totalheight=4cm]{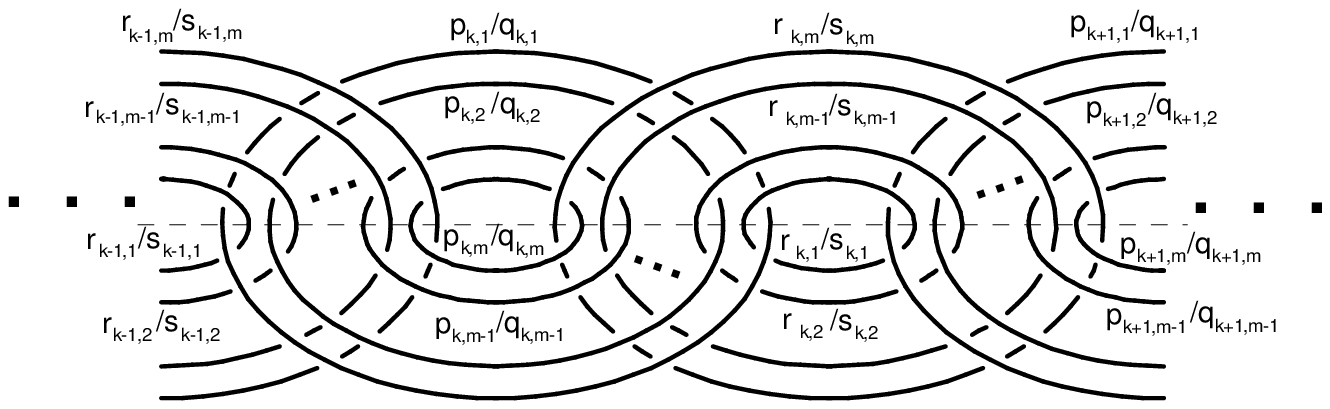}
 \end{center}
 \caption{Surgery presentation for generalized Takahashi manifolds.}
 \label{Fig. 1}
\end{figure}

In this section we define a family of manifolds which generalizes
Takahashi manifolds. For any pair of positive integers $m$ and
$n$, we consider the link ${\cal L}_{n,m} \subset \S^3$ with $2mn$
components presented in Figure \ref{Fig. 1}. All its components
$c_{i,j}$, $1\le i \le  2n$, $1 \le j \le m$, are unknotted
circles and they form $2n$ subfamilies of $m$ unlinked circles
$c_{i,j}$, $1 \le j \le m$, with a common center. We observe that
${\cal L}_{n,1}$ is the link ${\cal L}_n$ discussed above. The
link ${\cal L}_{n,m}$ has a cyclic symmetry of order $n$ which
permutes these $2n$ subfamilies of circles.

Consider the manifold obtained by Dehn surgery on $\S^3$, along
the link ${\cal L}_{n,m}$, such that the surgery coefficients $p_{k,j} /
q_{k,j}$ correspond to the components $c_{2k-1,j}$ and $r_{k,j} /
s_{k,j}$ correspond to the components $c_{2k,j}$, where $1 \le k
\le n$ and $1 \le j \le m$ (see Figure \ref{Fig. 1}). Without loss
of generality, we can always suppose $\gcd(p_{k,j},q_{k,j})=1$,
$\gcd(r_{k,j},s_{k,j})=1$ and $p_{k,j},r_{k,j}\ge 0$.

We will denote the resulting 3-manifold by $T_{n,m} (p_{k,j} /
q_{k,j} ; r_{k,j} / s_{k,j} )$. This manifold will be referred to as
a {\it generalized Takahashi manifold}, since for $m=1$ we get the
Takahashi manifolds introduced in \cite{Ta}.

The following theorem generalizes the result obtained in \cite{Ta}
for Takahashi manifolds.

\begin{theorem} \label{fundamental}
The fundamental group of the generalized Takahashi manifold
$T_{n,m} (p_{k,j} / q_{k,j} ; r_{k,j} / s_{k,j} )$ has the
following balanced presentation with $2nm$ ge\-ne\-rators
${A} = \{ a_{i,j} \}_{\, 1 \le i \le 2n \, , \,1 \le j \le m}$
and $2nm$ relations:
\vbox{
\begin{eqnarray}
\langle \, {A}  \mid
& a_{2k-1,j}^{-p_{k,j}} \, = \, a_{2k-2,j}^{s_{k-1,j}} \,
a_{2k-2,j+1}^{s_{k-1,j+1}} \, \cdots \, a_{2k-2,m}^{s_{k-1,m}} \,
a_{2k,m}^{-s_{k,m}} \, \cdots
\, a_{2k,j+1}^{-s_{k,j+1}} \, a_{2k,j}^{-s_{k,j}},  \nonumber \\
& a_{2k,j}^{-r_{k,j}} \, =
\, a_{2k+1,j}^{q_{k+1,j}} \,  a_{2k+1,j-1}^{q_{k+1,j-1}} \,
\cdots \, a_{2k+1,1}^{q_{k+1,1}} \, a_{2k-1,1}^{-q_{k,1}} \,
\cdots \, a_{2k-1,j-1}^{-q_{k,j-1}}  \, a_{2k-1,j}^{-q_{k,j}};
\nonumber \\ &  \qquad \qquad \qquad \qquad \qquad
 1 \le k \le n \, , \, 1 \le j \le m \, \rangle. \nonumber
\end{eqnarray}
}
\end{theorem}

\begin{proof}
Let $ X = \{ x_{k,j} \}_{\, 1\le k\le n, \, 1\le j\le m}$ and
$ Y = \{ y_{k,j} \}_{\, 1\le k\le n, \, 1\le j\le m}$ be sets of
Wirtinger generators of $\p_1(\S^3 \setminus \L_{n,m})$, according
to Figure \ref{Fig. 2}.

\bigskip\bigskip

\begin{figure}[bht]
 \begin{center}
 \includegraphics*[totalheight=4cm]{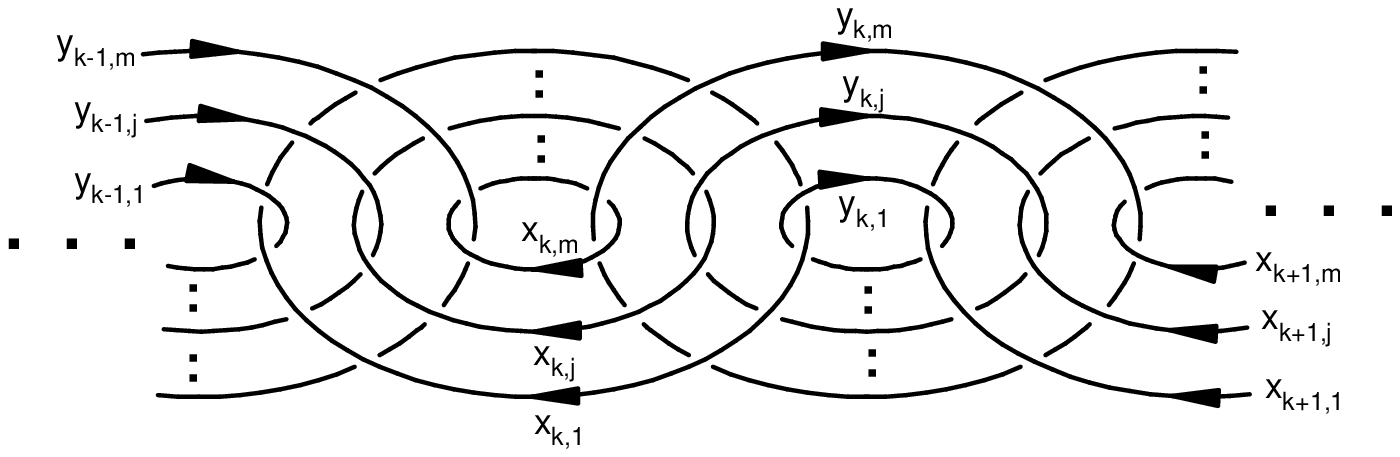}
 \end{center}
 \caption{Generators of $\p_1(\S^3 \setminus \L_{n,m})$.}
 \label{Fig. 2}
\end{figure}

Applying the Wirtinger algorithm we get the following presentation
for $\p_1(\S^3 \setminus \L_{n,m})$:
\begin{eqnarray}
\langle  X \cup  Y
\mid & y_{k,j}\cdots y_{k,m}y_{k-1,m}^{-1}\cdots
y_{k-1,j}^{-1}x_{k,j}y_{k-1,j}\cdots y_{k-1,m}y_{k,m}^{-1}\cdots
y_{k,j}^{-1}=x_{k,j} , \qquad\nonumber \\ & x_{k,j}\cdots
x_{k,1}x_{k+1,1}^{-1}\cdots x_{k+1,j}^{-1}y_{k,j}x_{k+1,j}\cdots
x_{k+1,1}x_{k,1}^{-1}\cdots x_{k,j}^{-1}=y_{k,j}  ; \qquad
\nonumber \\ & \qquad \qquad \qquad \qquad \qquad \qquad \qquad
\qquad 1\le k\le n\,,\,1\le j\le m \, \rangle. \nonumber
\end{eqnarray}
For every $k=1,\ldots,n$ and $j=1,\ldots,m$, let $h_{k,j}$ and
$l_{k,j}$ be the longitudes associated to the components of
$\L_{n,m}$ corresponding to the meridians $x_{k,j}$ and $y_{k,j}$
respectively (as usual we consider longitudes which are
homologically trivial in the complement of the relative
component). Then we have the relations: $$ h_{k,j} \, = \,
y_{k-1,j} \, y_{k-1,j+1} \, \cdots \, y_{k-1,m} \, y_{k,m}^{-1} \,
\cdots \, y_{k,j+1}^{-1} \, y_{k,j}^{-1}; $$ and $$ l_{k,j} \, =
\, x_{k+1,j} \, x_{k+1,j-1} \, \cdots \, x_{k+1,1} \, x_{k,1}^{-1}
\, \cdots \, x_{k,j-1}^{-1} \, x_{k,j}^{-1}. $$ Introducing
$H=\{h_{k,j}\}_{\, 1\le k\le n, \, 1\le j\le m}$ and
$L=\{l_{k,j}\}_{\, 1\le k\le n, \, 1\le j\le m}$, we obtain
the following new presentation for $\p_1(\S^3 \setminus
\L_{n,m})$:
\begin{eqnarray}
\langle X \cup Y \cup H \cup L \mid & [x_{k,j},h_{k,j}]=1,  \quad
h_{k,j}=y_{k-1,j}\cdots y_{k-1,m}y_{k,m}^{-1}\cdots y_{k,j}^{-1}
\,, \nonumber \\ & [y_{k,j},l_{k,j}]=1, \quad
l_{k,j}=x_{k+1,j}\cdots x_{k+1,1}x_{k,1}^{-1}\cdots
x_{k,j}^{-1}\,;\, \nonumber \\ & \qquad \qquad \qquad \qquad
\qquad \qquad {1\le k\le n\,,\,1\le j\le m} \rangle. \nonumber
\end{eqnarray}
Therefore, the fundamental group of $T_{n,m} (p_{i,j} / q_{i,j} ; r_{i,j} /
s_{i,j} )$ admits the presentation:
\begin{eqnarray}
\langle X \cup Y \cup H \cup L \mid & [x_{k,j}, h_{k,j}] = 1,
\quad h_{k,j} = y_{k-1,j}\cdots y_{k-1,m}y_{k,m}^{-1}\cdots
y_{k,j}^{-1}\,,\, \nonumber \\ & [y_{k,j},l_{k,j}] = 1, \quad
l_{k,j} = x_{k+1,j}\cdots x_{k+1,1}x_{k,1}^{-1}\cdots
x_{k,j}^{-1}\,;\, \nonumber \\ & x_{k,j}^{p_{k,j}}
h_{k,j}^{q_{k,j}} = 1, \quad y_{k,j}^{r_{k,j}} l_{k,j}^{s_{k,j}} =
1, \quad {1\le k\le n\,,\,1\le j\le m} \rangle . \nonumber
\end{eqnarray}

Since $\gcd(p_{k,j},q_{k,j})=1$ and $\gcd(r_{k,j},s_{k,j})=1$,
there exist certain integers $u_{k,j}$, $v_{k,j}$, $w_{k,j}$ and
$z_{k,j}$ such that $q_{k,j} u_{k,j} - p_{k,j} v_{k,j} = 1$ and
$s_{k,j} w_{k,j} - r_{k,j} z_{k,j} = 1$.

For $k=1,\ldots,n$ and $j=1,\ldots,m$ we define $$ a_{2k-1,j} =
x_{k,j}^{u_{k,j}} h_{k,j}^{v_{k,j}}, \qquad a_{2k,j} =
y_{k,j}^{w_{k,j}} l_{k,j}^{z_{k,j}}. $$ Since $x_{k,j}$ and
$h_{k,j}$ (resp. $y_{k,j}$ and $h_{k,j}$) commute, we have
\begin{eqnarray}
a_{2k-1,j}^{q_{k,j}} = x_{k,j}(x_{k,j}^{p_{k,j} v_{k,j}}
h_{k,j}^{q_{k,j} v_{k,j}}) = x_{k,j}, \nonumber \\
a_{2k-1,j}^{-p_{k,j}} = (x_{k,j}^{-p_{k,j} u_{k,j}}
h_{k,j}^{-q_{k,j} u_{k,j}}) h_{k,j} = h_{k,j}, \nonumber \\
a_{2k,j}^{s_{k,j}} = y_{k,j}(y_{k,j}^{r_{k,j} z_{k,j}}
h_{k,j}^{s_{k,j} z_{k,j}}) = y_{k,j}, \nonumber \\
a_{2k,j}^{-r_{k,j}} = ( y_{k,j}^{-r_{k,j} w_{k,j}}
l_{k,j}^{-s_{k,j} w_{k,j}}) l_{k,j} = l_{k,j}. \nonumber
\end{eqnarray}
Using these relations we can eliminate all the generators of the previous
presentation of $T_{n,m} (p_{k,j} / q_{k,j} ; r_{k,j} / s_{k,j}
)$, replacing them with the set $\{ a_{i,j}\}_{\, 1 \le i \le 2n,
\, 1 \le j \le m}$. The first four types of relations of the above
presentation disappear and the statement is obtained.
\end{proof}

\medskip

When the surgery coefficients are $n$-periodic, i.e.
$p_{k,j}=p_j$, $q_{k,j}=q_j$, $r_{k,j}=r_j$, and $s_{k,j}=s_j$,
the resulting manifold $T_{n,m} (p_j / q_j ; r_j / s_j )$ is said
to be a {\it generalized periodic ($n$-periodic) Takahashi
manifold}.

\begin{corollary} \label{fundamental periodic}
The fundamental group of the generalized periodic Takahashi
manifold $T_{n,m} (p_j / q_j ; r_j / s_j )$ admits the presentation
\begin{eqnarray}
\langle \, \{ a_{i,j} \}_{\, 1 \le i \le 2n \, , \,1 \le j \le m}
\mid & a_{2k-1,j}^{-p_j} \, = \, a_{2k-2,j}^{s_j} \, \cdots \,
a_{2k-2,m}^{s_m} \, a_{2k,m}^{-s_m} \, \cdots \, a_{2k,j}^{-s_j},
\nonumber \\ & a_{2k,j}^{-r_j} \, = \, a_{2k+1,j}^{q_j} \, \cdots
\, a_{2k+1,1}^{q_1} \, a_{2k-1,1}^{-q_1} \, \cdots \,
a_{2k-1,j}^{-q_j};  \nonumber \\ &  \qquad \qquad \qquad \qquad
\qquad
 1 \le k \le n \, , \, 1 \le j \le m \, \rangle. \nonumber
\end{eqnarray}
\end{corollary}

\section{Covering properties of generalized Takahashi manifolds}

\begin{figure}[bht]
 \begin{center}
 \includegraphics*[totalheight=6cm]{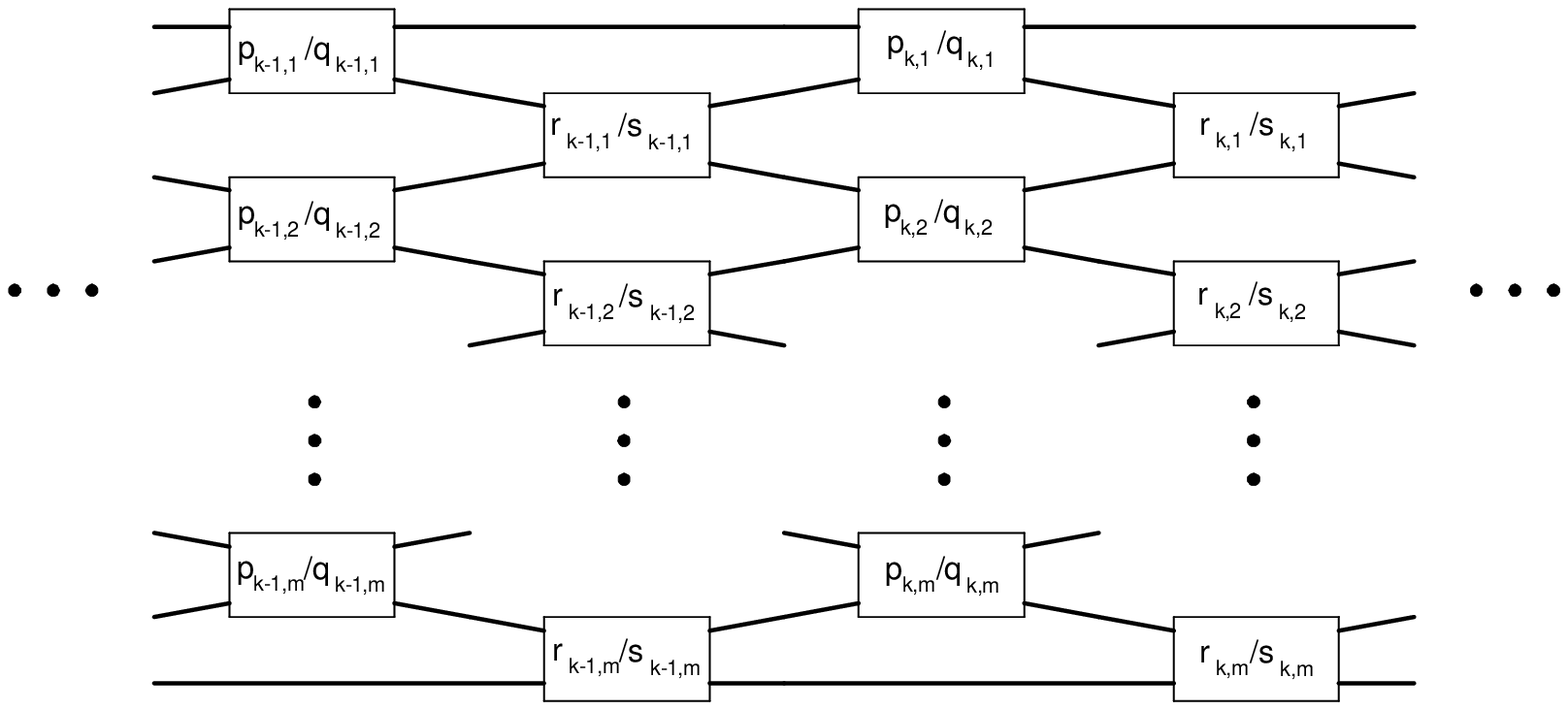}
 \end{center}
 \caption{The link ${\cal K}_{n,m}(p_{k,j}/q_{k,j};
r_{k,j}/s_{k,j})$.}
 \label{Fig. 3}
\end{figure}

We will define a new family of links in $\S^3$. For any pair of
integers $n,m>0$ consider two pairs of coprime integers
$(p_{k,j},q_{k,j})$ and $(r_{k,j},s_{k,j})$, where $k=1, \ldots ,
n$ and  $j=1, \ldots , m$. Let ${\cal K}_{n,m}(p_{k,j}/q_{k,j};
r_{k,j}/s_{k,j})$ be the closure of the rational braid on $2m+1$
strings with rational tangles \cite{BZ} $p_{k,j}/q_{k,j}$ and
$r_{k,j}/s_{k,j}$ indicated in Figure \ref{Fig. 3}.

As a generalization of the results from \cite{KV,RS,Ta}, we get:

\begin{theorem} \label{twofold}
The generalized Takahashi manifold $T_{n,m} (p_{k,j} / q_{k,j} ;
r_{k,j} / s_{k,j} )$ is the 2-fold covering of $\S^3$, branched
over the link ${\cal K}_{n,m} (p_{k,j} / q_{k,j} ; r_{k,j} /
s_{k,j} )$.
\end{theorem}

\begin{proof} From Figure \ref{Fig. 1} we see that the link $\L_{n,m}$
admits a strongly invertible involution $\t$ whose axis (pictured
with dashed line) intersects each component of the link
in two points. Thus, in virtue of the Montesinos theorem \cite{Mo},
the manifold $T_{n,m} (p_{k,j} / q_{k,j} ; r_{k,j} / s_{k,j} )$
can be obtained as the $2$-fold covering of $\S^3$, branched over some
link.

\begin{figure}[bht]
 \begin{center}
 \includegraphics*[totalheight=6.0cm]{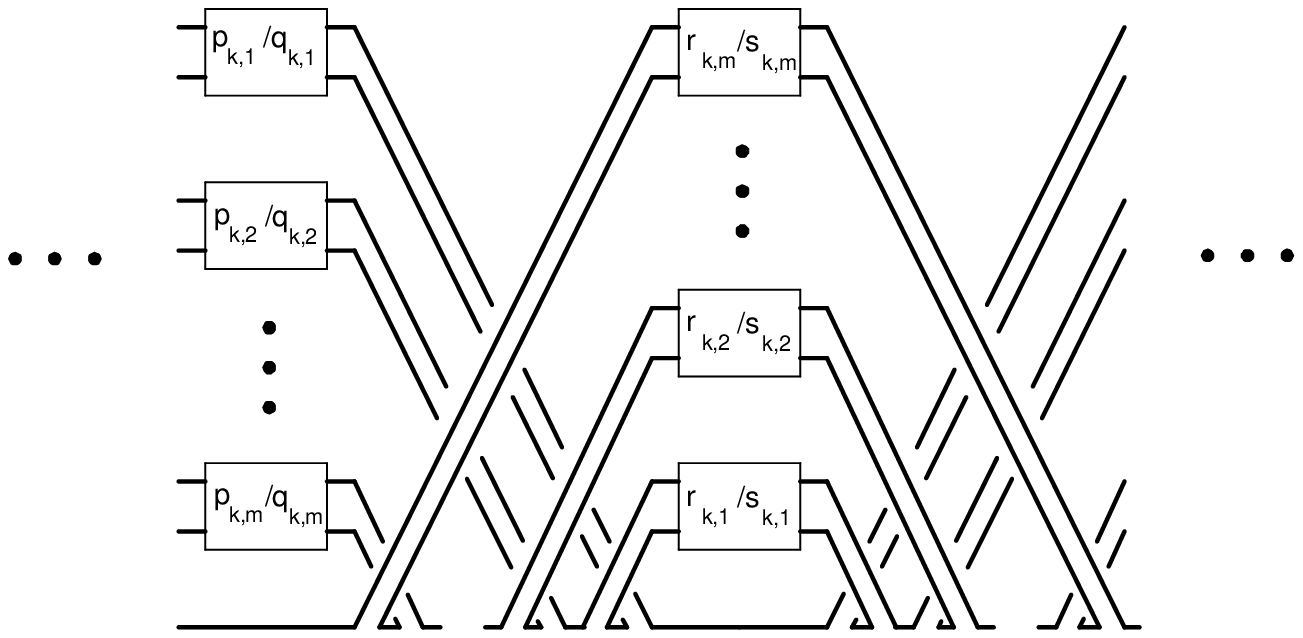}
 \end{center}
 \caption{}
 \label{Fig. 4}
\end{figure}

Applying the Montesinos algorithm, we get the link depicted in
Figure \ref{Fig. 4}. Obviously, this branching set is equivalent
to the link presented in Figure \ref{Fig. 3}. \end{proof}

\medskip

In particular, if the surgery coefficients are $n$-periodic, i.e.
$p_{k,j}=p_j$, $q_{k,j}=q_j$, $r_{k,j}=r_j$, and $s_{k,j}=s_j$,
the link ${\cal K}_{n,m} (p_{j} / q_{j} ; r_{j} / s_{j} )$ is also
$n$-periodic. Note that ${\cal K}_{n,1}(1;-1)$ is an alternating
link with $2n$ double-crossings, which is the closure of a 3-string
braid, referred to as a Turk head link in \cite{MV}: ${\cal
K}_{2,1}(1;-1)$ is the figure-eight knot and ${\cal
K}_{3,1}(1;-1)$ are the Borromean rings.

\begin{corollary} \label{twofold periodic}
The generalized periodic Takahashi manifold $T_{n,m} (p_{j} /
q_{j} ; r_{j} / s_{j} )$ is the 2-fold covering of $\S^3$,
branched over the periodic link ${\cal K}_{n,m} (p_j / q_j ; r_j /
s_j )$.
\end{corollary}

\begin{figure}[bht]
 \begin{center}
 \includegraphics*[totalheight=10cm]{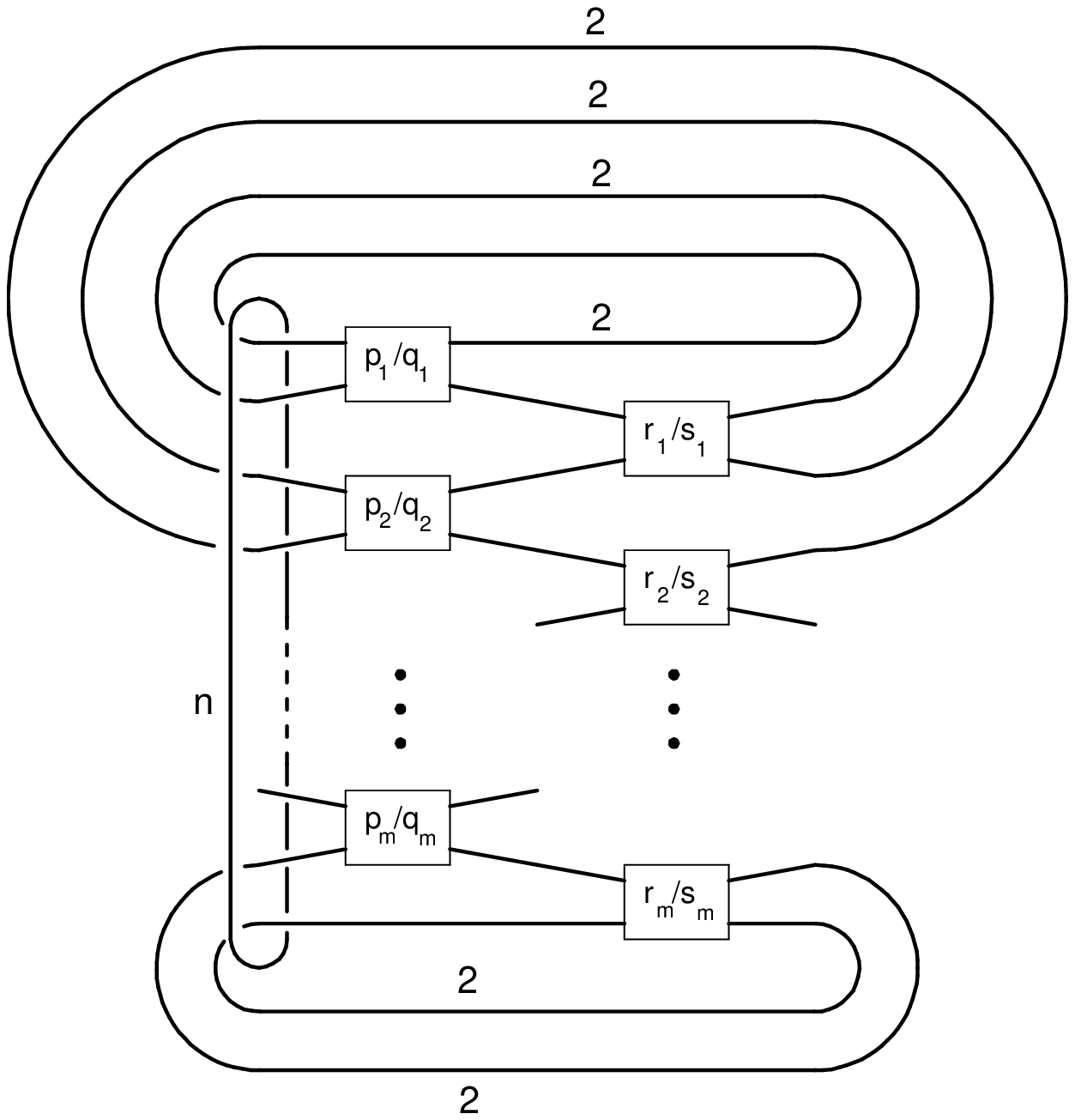}
 \end{center}
 \caption{The link ${\cal Q}_{n,m} (p_j / q_j ; r_j / s_j )$.}
 \label{Fig. 6}
\end{figure}

In other words, $T_{n,m} (p_{j} / q_{j} ; r_{j} / s_{j} )$ is the
$\Z_2$-covering of the orbifold $\S^3 ({\cal K}_{n,m} (p_j / q_j ;
r_j / s_j ))$ whose underlying space is $\S^3$ and whose singular set is
${\cal K}_{n,m} (p_j / q_j ; r_j / s_j )$, with singularity indices
$2$. Since the singular set of the orbifold is $n$-periodic, there
is a natural action of a cyclic group $\Z_n$ such that the
quotient orbifold is $\S^3 ({\cal Q}_{n,m} (p_j / q_j ; r_j / s_j
))$, where the singular set is the link pictured in Figure
\ref{Fig. 6} and the indices of singularity are: $2$ on the
components which are images of ${\cal K}_{n,m} (p_j / q_j ; r_j /
s_j )$ and $n$ on the unknotted component. Note that the part of
the singular set having index $2$ can be obtained as a connected
sum of $2m$ two-bridge links corresponding to the rational tangles
$p_1/q_1, r_1/s_1, \ldots , p_m/q_m, r_m/s_m$.

Therefore we get the following statement.

\begin{corollary} \label{twofold cyclic}
The generalized periodic Takahashi manifold $T_{n,m} (p_{j} /
q_{j} ; r_{j} / s_{j} )$ is the $\Z_2 \oplus \Z_n$-covering of the
orbifold $\S^3 ({\cal Q}_{n,m} (p_j / q_j ; r_j / s_j ))$.
\end{corollary}

The following theorem extends to generalized periodic Takahashi
manifolds the result given in \cite{Mu2} for periodic Takahashi
manifolds.

\begin{theorem} \label{cyclic}
The generalized periodic Takahashi manifold $T_{n,m} (p_j / q_j ;
r_j / s_j )$ is the $n$-fold cyclic covering of the
connected sum of $2m$ lens spaces $L(p_1,q_1) \# L(r_1,s_1) \#
\cdots \# L(p_m,q_m) \# L(r_m,s_m)$, branched over a knot which
does not depend on $n$.
\end{theorem}

\begin{proof}
Both the link $\L_{n,m}$ and the surgery coefficients defining the
ma\-ni\-fold $T_{n,m} (p_j / q_j ; r_j / s_j )$ (and so, also the manifold)
are invariant with respect to an obvious rotation symmetry $\rho$
of order $n$. Denote by $\langle \rho \rangle$ the cyclic group of
order $n$ generated by this rotation. Observe that the fixed-point
set of the action of $\langle \rho \rangle$ on $\S^3$ is a trivial
knot disjoint from $\L_{n,m}$. Therefore, we have an action of
$\langle \rho \rangle$ on $T_{n,m} (p_j / q_j ; r_j / s_j )$, with
a knot $K=K(\rho)$ as fixed-point set. The underlying space of the
quotient orbifold $T_{n,m} (p_j / q_j ; r_j / s_j ) / \langle \rho
\rangle$ is precisely the manifold $T_{1,m} (p_j / q_j ; r_j / s_j
)$, which can be obtained by Dehn surgery on $\S^3$, with
coefficients $p_j / q_j$ and $r_j / s_j$, $j=1, \ldots , m$, along
the $2m$-component link $\L_{1,m}$ depicted in Figure \ref{Fig.
5}.

\begin{figure}[bht]
 \begin{center}
 \includegraphics*[totalheight=7cm]{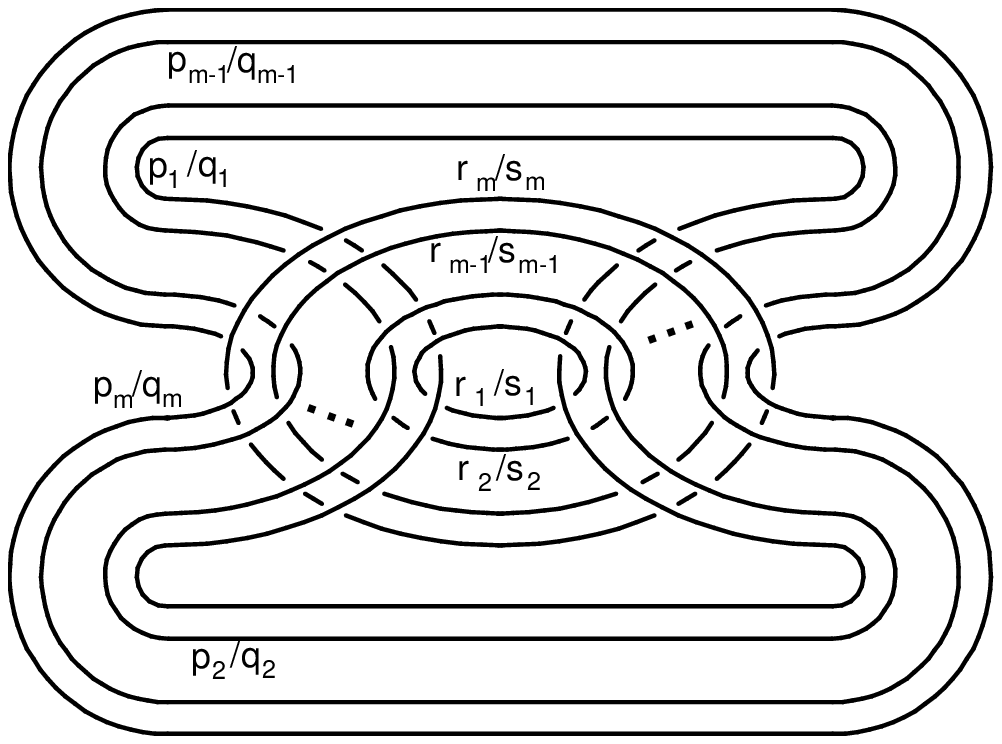}
 \end{center}
 \caption{The link $\L_{1,m}$.}
 \label{Fig. 5}
\end{figure}

The components of $\L_{1,m}$ are unlinked, unknotted, and form a
trivial link with $2m$ components. Therefore the underlying space of the quotient
or\-bi\-fold is homeomorphic to the connected sum of $2m$ lens
spaces $L(p_1,q_1) \# L(r_1,q_1) \# \cdots \# L(p_m,q_m) \#
L(r_m,s_m)$
 (see \cite[p.~260]{Ro}).
Moreover, it is obvious from the action of $\rho$ that the
singular set $K$ of the quotient orbifold is a knot which does not
depend on $n$.
\end{proof}

\medskip

Denote by ${\cal O}_{n,m} (p_j / q_j ; r_j / s_j )= T_{n,m} (p_j /
q_j ; r_j / s_j ) / \langle \rho \rangle$ the orbifold from the
proof of Theorem \ref{cyclic}, whose underlying space is the
connected sum of $2m$ lens spaces
$L(p_1,q_1)\#L(r_1,s_1)\#\dots\#L(p_m,q_m)\#L(r_m,s_m)$.

\begin{corollary} \label{commutative}
The following commutative diagram holds for each generalized
periodic Takahashi manifold.

 \begin{center}
  \unitlength=0.5mm
  \begin{picture}(90,75)(0,10)
  \put(40,80){\makebox(0,0)[cc]{$T_{n,m} (p_j / q_j ; r_j / s_j )$}}
  \put(50,70){\vector(2,-1){30}}
  \put(30,70){\vector(-2,-1){30}}
  \put(68,68){\makebox(0,0)[cc]{$n$}} 
  \put(7,68){\makebox(0,0)[cc]{$2$}} 
  \put(0,45){\makebox(0,0)[cc]{$\S^3({\cal K}_{n,m} (p_j / q_j ; r_j / s_j ))$}}
  \put(80,45){\makebox(0,0)[cc]{${\cal O}_{n,m} (p_j / q_j ; r_j / s_j )$}}
  \put(0,35){\vector(2,-1){30}}
  \put(80,35){\vector(-2,-1){30}}
  \put(7,25){\makebox(0,0)[cc]{$n$}} 
  \put(68,25){\makebox(0,0)[cc]{$2$}} 
  \put(40,10){\makebox(0,0)[cc]{$\S^3({\cal Q}_{n,m} (p_j / q_j ; r_j / s_j ))$}}
  \end{picture}
  \end{center}

\end{corollary}

\begin{proof} From Figure \ref{Fig. 1} we see that ${\cal L}_{n,m}$ admits an
invertible involution $\tau$ whose axis intersects each component
in two points and the rotation symmetry $\rho$ of order $n$ which
was discussed in Theorem \ref{cyclic}. These symmetries induce
symmetries (also denoted by $\tau$ and $\rho$) of the generalized
periodic Takahashi manifold $T = T_{n,m} (p_j / q_j ; r_j / s_j )$, such that
$\langle \tau ,
\rho \rangle \cong  \langle \tau \rangle
\oplus \langle \rho \rangle \cong \Z_2 \oplus \Z_n $. As mentioned above, $\rho$ induces
the symmetry (also denoted by $\rho$) of the orbifold $T /
\langle \tau \rangle$ (whose singular set is given by Corollary
\ref{twofold periodic}), and the covering $T \to (T / \langle \tau
\rangle )/ \langle \rho \rangle$ is given by Corollary \ref{twofold
cyclic}. The covering $T \to T / \langle \rho \rangle$ is given by
Theorem \ref{cyclic}. As we see from Figure \ref{Fig. 5}, $\tau$
induces the strongly invertible involution (also denoted by
$\tau$) of the link ${\cal L}_{1,m}$. Using the Montesinos
algorithm we see that $(T / \langle \rho \rangle ) / \langle \tau
\rangle = \S^3({\cal Q}_{n,m} (p_j / q_j ; r_j / s_j ))$ (note
that the part of the singular set of $\S^3 ({\cal Q}_{n,m} (p_j /
q_j ; r_j / s_j ))$ having index $2$ can be obtained as a
connected sum of $2m$ two-bridge links corresponding to the rational
tangles $p_1/q_1, r_1/s_1, \ldots ,p_m/q_m, r_m/s_m$).
\end{proof}

\medskip

\section{Cyclic branched coverings of 2-bridge knots}

In this section we show that generalized periodic Takahashi
manifolds contain the whole class of cyclic branched coverings of
two-bridge knots.
In the following we use the Conway notation for two-bridge knots
(see \cite{Co}).

\begin{theorem} \label{twobridge}
The generalized periodic Takahashi manifold $T_{n,m}(1/q_j;
1/s_j)$ is the $n$-fold cyclic branched covering of the two-bridge
knot corresponding to the Conway parameters
$[-2q_1,2s_1,\ldots,-2q_m,2s_m]$.
\end{theorem}

\begin{proof} From Theorem \ref{cyclic}, $T_{n,m}(1/q_j; 1/s_j)$ is the $n$-fold
cyclic co\-ve\-ring of $\S^3$, branched over a knot $K$. Figures
\ref{Fig. 7}--\ref{Fig. 13} shows how to deform $K$ to a
Conway's normal form of a two-bridge knot with Conway parameters
$[-2q_1,2s_1,\ldots,-2q_m,2s_m]$ by ambient isotopy (from Figure
\ref{Fig. 7} to Figure \ref{Fig. 10}) and surgery calculus
\cite{Ro} (from Figure \ref{Fig. 11} to Figure \ref{Fig. 13}).
\end{proof}

\medskip

\noindent {\bf Remark} As a consequence of Theorem
\ref{twobridge}, the generalized periodic Takahashi manifold
$T_{n,m}(1/q_j; 1/s_j)$ is homeomorphic to the Lins-Mandel
manifold $S(n,a,b,1)$ \cite{LM,Mu1}, the Minkus manifold
$M_n(a,b)$ \cite{Mi} and the Dunwoody manifold
$M((a-1)/2,0,1,b/2,n,-q_{\s})$ \cite{Du,GM}, where $a/b$ is
defined by (1).

\medskip

Because every 2-bridge knot admits a Conway representation with an
even number of even parameters (see,  Exercise 2.1.14 of
\cite{Ka}), we have the following property.

\begin{corollary} \label{all twobridge}
The family of generalized periodic Takahashi manifolds contains
all cyclic branched coverings of two-bridge knots.
\end{corollary}

\begin{figure}[bht]
 \begin{center}
 \includegraphics*[totalheight=4cm]{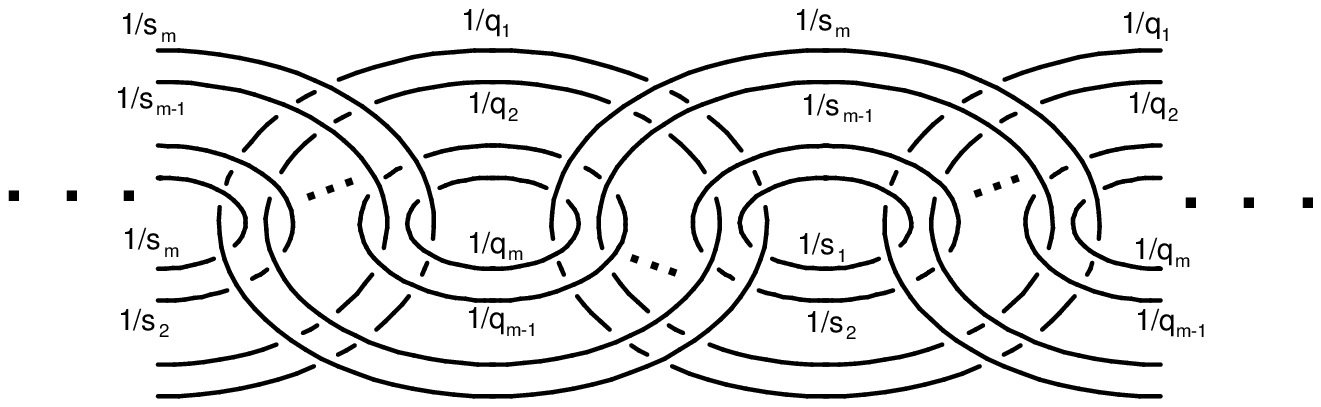}
 \end{center}
 \caption{Surgery presentation for
 $\wh C_n\big(-2q_1+\frac{1}{2s_1+\cdots + \frac{1}{-2q_m+\frac{1}{2s_m}}}\big)$.}
 \label{Fig. 00}
\end{figure}

From Theorem \ref{twobridge} we can easily get the surgery
presentation for the $n$-fold cyclic branched covering $\wh C_n(a/b)$ of
the two-bridge knot, with Conway parameters
$[-2q_1,2s_1,\ldots,-2q_m,2s_m]$, depicted in Figure \ref{Fig. 00}.

\section{Cyclically presented fundamental groups}

A cyclic presentation for the fundamental groups of cyclic
branched coverings of two-bridge knots is obtained by J. Minkus
(see Theorem 10 of \cite{Mi}). Corollary \ref{fundamental
periodic} and Theorem \ref{twobridge} allow us to obtain a
different cyclic presentation for such groups. Note that
explicit cyclic presentations different from the above are listed
in the Appendix of \cite{CRS}, for two-bridge knots up to nine
crossings.

\begin{theorem} \label{cyclicpresentation}
Let $\wh C_n(a/b)$ be the $n$-fold cyclic branched covering of the
two-bridge knot ${\bf b}(a/b)$, with $a/b$ given by formula (1).
Then its fundamental group has the following cyclic presentation:
$$ \pi_1 (\wh C_n (a/b)) \, = \, \langle x_1, \dots , x_n \, | \,
w_{a/b} (x_i, \dots , x_{i+n-1}) = 1, \qquad i=1, \dots , n
\rangle , $$ where
$$
w_{a/b} (x_i, \ldots , x_{i+n-1} = b^{-s_m}_{i+1,m} d_{i+1}
b^{s_m}_{1,m}
$$
for $i=1, \ldots, n$ (indices mod $n$). The right parts of these
formulas are defined by the recurrent rule
$$ d_{k,j} = b_{k,j-1}^{-s_{j-1}} d_{k,j-1} b_{k-1,j-1}^{s_{j-1}},
\qquad b_{k,j} = d_{k,j}^{q_j} b_{k,j-1} d_{k+1,j}^{-q_j}, \qquad
j=2, \dots , m $$ and $$ b_{k,1} = d_{k,1}^{q_1} d_{k+1,1}^{-q_1},
$$ where $x_k = d_{k,1}$, for $k=1, \dots , n$.
\end{theorem}

\begin{proof} From Corollary \ref{fundamental periodic} and
Theorem \ref{twobridge}, the group $\pi_1(\wh C_n(a/b))$ is
generated by the $2nm$ elements $\{ a_{i,j} \}_{i=1, \dots,2n, \,
j=1,\dots,m}$ and has relations of two types:
\begin{eqnarray}
a_{2k-1,j}^{-1}  = &
\left( a_{2k-2,j}^{s_j} a_{2k-2,j+1}^{s_{j+1}} \cdots
a_{2k-2,m}^{s_m} \right)
\cdot
\left( a_{2k,j}^{s_j} a_{2k,j+1}^{s_{j+1}} \cdots
a_{2k,m}^{s_m} \right)^{-1}, \nonumber \\
a_{2k,j}^{-1} = &
\left( a_{2k+1,j}^{q_j} a_{2k+1,j-1}^{q_{j-1}} \cdots
a_{2k+1,1}^{q_1} \right)
\cdot
\left( a_{2k-1,j}^{q_j} a_{2k-1,j-1}^{q_{j-1}} \cdots
a_{2k-1,1}^{q_1} \right)^{-1}, \nonumber
\end{eqnarray}
where $k=1, \dots ,n$ and $j=1, \dots , m$, and all the indices
are taken mod $2n$ and $m$ respectively. Denote
$b_{k,j}=a_{2k,j}$ and $d_{k,j}=a_{2k-1,j}$ for $k=1, \dots ,n$
and $j=1, \dots , m$. Then we have $2nm$ relations of the two
following types: $$ d_{k,j} = b_{k,j}^{s_j} b_{k,j+1}^{s_{j+1}}
\cdots b_{k,m}^{s_m} b_{k-1,m}^{-s_m} \cdots
b_{k-1,j+1}^{-s_{j+1}} b_{k-1,j}^{-s_j} $$ and $$ b_{k,j} =
d_{k,j}^{q_j} d_{k,j-1}^{q_{j-1}} \cdots d_{k,1}^{q_1}
d_{k+1,1}^{-q_1} \cdots d_{k+1,j-1}^{-q_{j-1}} b_{k+1,j}^{-q_j}
.$$ Therefore, the defining relations for the group are: $$
b_{k,m}^{-s_m} d_{k,m} b_{k-1,m}^{s_m} = 1, \qquad d_{k,j+1} =
b_{k,j}^{-s_j} d_{k,j} b_{k-1,j}^{s_j}, \quad j=1, \dots, m-1 ,$$
and $$ b_{k,1} = d_{k,1}^{q_1} d_{k+1,1}^{-q_1}, \qquad b_{k,j} =
d_{k,j}^{q_j} b_{k,j-1} d_{k+1,j}^{-q_j}, \quad j=2, \dots, m ,$$
for $k=1, \dots ,n$. Denoting $x_k=d_{k,1}$, $k=1, \dots ,n$, we
will eliminate all other generators in the following order:
$b_{k,1}, d_{k,2}, b_{k,2}, \dots , d_{k,m}, b_{k,m}$ according to
the above formulae. At the end of this process we will get $n$
relations arising from $b_{k,m}^{-s_m} d_{k,m} b_{k-1,m}^{s_m} =
1$. That completes the proof.
\end{proof}

\medskip

We will illustrate the obtained result for the cases $m=1$ and $m=2$.

\medskip

If $m=1$, then $a/b = -2q + \frac{\displaystyle 1}{\displaystyle
2s}$, and $\wh C_n (a/b) = T_{n,1}(1/q,1/s)$. This case,
corresponding to a Takahashi manifold, was discussed in
\cite{KKV} and \cite{KV}. Using notations $b_k=b_{k,1}$ and
$d_k=d_{k,1}$ for $k=1, \dots ,n$, we get
\begin{eqnarray}
\pi_1 (T_{n,1}(1/q,1/s)) \, = \, & \nonumber \\
 \langle b_1, \dots , b_n, d_1, \dots ,d_n \, & \, |  \quad
 b_{k}^{-s} d_{k} b_{k}^{s} = 1, \qquad
 b_{k} = d_{k}^{q} d_{k+1}^{-q}, \qquad   k=1,\dots ,n \rangle . \nonumber
\end{eqnarray}
Hence
$$
\pi_1 (T_{n,1}(1/q,1/s)) \, = \,
\langle
x_1, \dots , x_n \, | \,
(x_k^q x_{k+1}^{-q})^{-s} x_k (x_{k-1}^q x_k^{-q})^s = 1,
 \quad k=1,\dots ,n
\rangle .
$$

For example, if $q=-1$ and $s=1$ then $a/b=5/2$, that corresponds
to the figure-eight knot $4_1$ \cite{BZ}. So, its $n$-fold cyclic branched
covering has the fundamental group with the cyclic
presentation
$$
\pi_1 (T_{n,1}(-1,1)) \, = \, \langle x_1, \dots , x_n \, | \, \,
x_{k+1}^{-1} x_k^2 x_{k-1}^{-1} x_k  = 1 ,  \quad k=1, \dots ,n
\rangle
$$
(compare with \cite{CRS, KKV, KV}).

For $m=2$ we get
$$
\pi_1 (T_{n,2}(1/q_1,1/q_2;1/s_1,1/s_2)) \,  = \qquad  \qquad
\qquad \qquad \qquad \qquad \qquad
$$
$$
\langle
b_{1,1}, \dots , b_{n,1}, b_{1,2}, \dots, b_{n,2} ,
d_{1,1}, \dots , d_{n,1}, d_{1,2}, \dots, d_{n,2}
| \qquad \qquad \qquad \qquad \qquad
$$
$$
b_{k,2}^{-s_2} d_{k,2} b_{k-1,2}^{s_2} = 1 , \qquad
d_{k,2} = b_{k,1}^{-s_1} d_{k,1} b_{k-1,1}^{s_1} , \qquad
b_{k,1} = d_{k,1}^{q_1} d_{k+1,1}^{-q_1} , \qquad
$$
$$
b_{k,2} = d_{k,2}^{q_2} b_{k,1} d_{k+1,2}^{-q_2}  ,
\qquad  \qquad  k=1, \dots ,n \rangle .
$$

Denote $x_k=d_{k,1}$, then $b_{k,1}=x_{k,1}^{q_1} x_{k+1}^{-q_1}$.
Therefore
$$
d_{k,2} = (x_k^{q_1} x_{k+1}^{-q_1})^{-s_1} x_k
(x_{k-1}^{q_1} x_{k}^{-q_1})^{s_1}
$$
and
$$
b_{k,2} =
\left[ (x_k^{q_1} x_{k+1}^{-q_1})^{-s_1} x_k
(x_{k-1}^{q_1} x_{k}^{-q_1})^{s_1} \right]^{q_2}
 (x_k^{q_1} x_{k+1}^{-q_1})
\left[ (x_{k+1}^{q_1} x_{k+2}^{-q_1})^{-s_1} x_{k+1}
(x_{k}^{q_1} x_{k+1}^{-q_1})^{s_1} \right]^{-q_2} .
$$

Define
$$
w_{a/b} (x_{k-2}, x_{k-1}, x_{k}, x_{k+1}, x_{k+2}) =
$$
$$
\left[ \left[ (x_k^{q_1} x_{k+1}^{-q_1})^{-s_1} x_k
(x_{k-1}^{q_1} x_{k}^{-q_1})^{s_1} \right]^{q_2}
 x_k^{q_1} x_{k+1}^{-q_1}
\left[ (x_{k+1}^{q_1} x_{k+2}^{-q_1})^{-s_1} x_{k+1}
(x_{k}^{q_1} x_{k+1}^{-q_1})^{s_1} \right]^{-q_2} \right]^{-s_2}
$$
$$
\cdot (x_k^{q_1} x_{k+1}^{-q_1})^{-s_1} x_k
(x_{k-1}^{q_1} x_{k}^{-q_1})^{s_1}
$$
$$
\cdot \left[ \left[ (x_{k-1}^{q_1} x_{k}^{-q_1})^{-s_1} x_{k-1}
(x_{k-2}^{q_1} x_{k-1}^{-q_1})^{s_1} \right]^{q_2}
 x_{k-1}^{q_1} x_{k}^{-q_1}
\left[ (x_{k}^{q_1} x_{k+1}^{-q_1})^{-s_1} x_{k}
(x_{k-1}^{q_1} x_{k}^{-q_1})^{s_1} \right]^{-q_2} \right]^{s_2}.
$$

Therefore, we get the following cyclic presentation for the
fundamental group of the $n$-fold cyclic branched covering of the
two-bridge knot ${\bf b}(a/b)$ corresponding to $ [-2q_1, 2s_1, -2q_2,
2s_2] $:
$$
\pi_1 (T_{n,2}(1/q_1,1/q_2;1/s_1,1/s_2)) \,  = \, $$ $$
\langle x_{1}, \dots , x_{n} \quad | \quad w_{a/b} (x_{k-2},
x_{k-1}, x_{k}, x_{k+1}, x_{k+1}) = 1, \qquad k=1, \dots , n
\rangle ,
$$
where all the indices are mod $n$.

For example, for $q_1=q_2=-1$ and
$s_1=s_2=1$ we get $a/b=29/12$, that corresponds to the knot
$8_{12}$. So, its $n$-fold cyclic branched covering has the
fundamental group with the following cyclic presentation:
\begin{eqnarray}
\langle x_1, \dots, x_n \, | \, &
x_{k+1}^{-1} x_{k} x_{k+1}^{-2} x_{k+2} x_{k+1}^{-1} x_{k} x_{k+1}^{-1}
x_{k}^{2} x_{k-1}^{-1} x_{k} x_{k+1}^{-1} &  \nonumber \\
& \cdot
x_{k}^{2} x_{k-1}^{-1} x_{k} x_{k-1}^{-1} x_{k-2} x_{k-1}^{-2}
x_{k} x_{k-1}^{-1} x_{k} x_{k+1}^{-1} x_k^2 x_{k-1}^{-1} x_k =1,
&  k=1, \dots ,n \rangle . \nonumber
\end{eqnarray}





\begin{figure}
\begin{center}
\includegraphics*[totalheight=8cm]{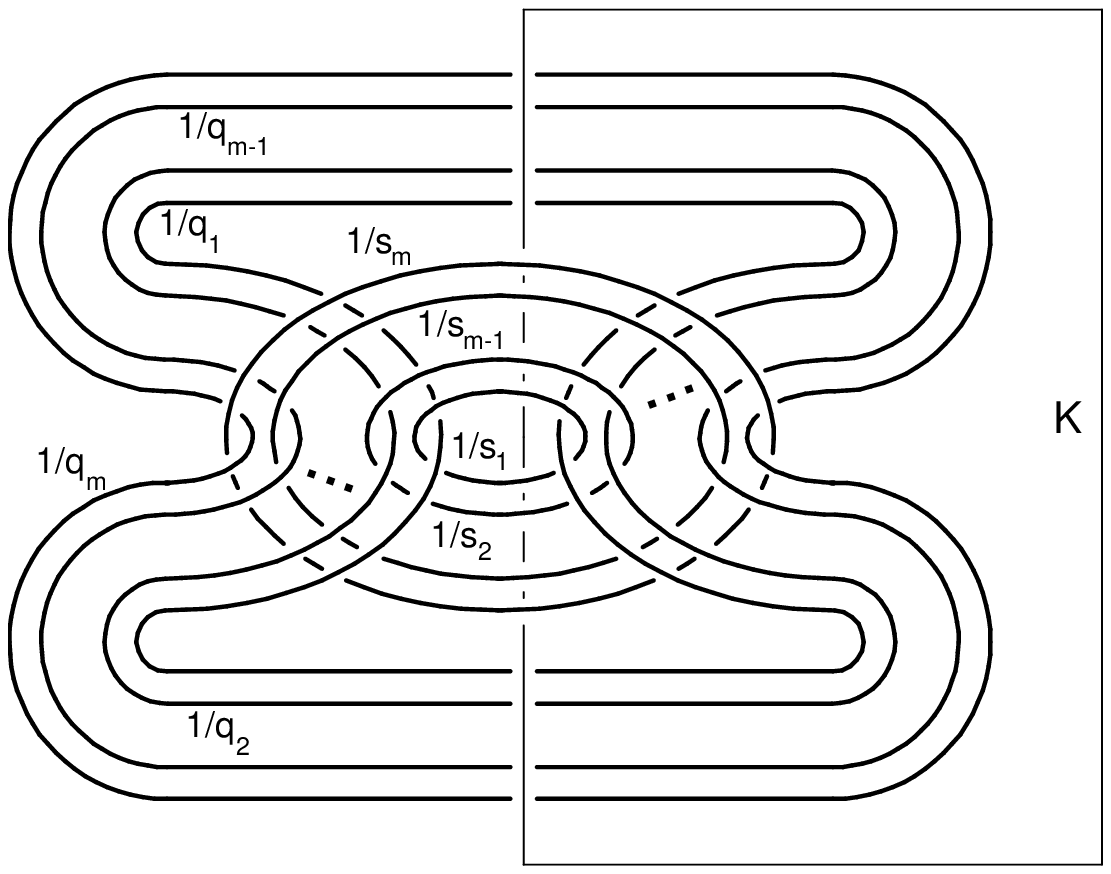}
\end{center}
\caption{$\L_{1,m}\cup K$.} \label{Fig. 7}
\end{figure}

\begin{figure}
 \begin{center}
 \includegraphics*[totalheight=8cm]{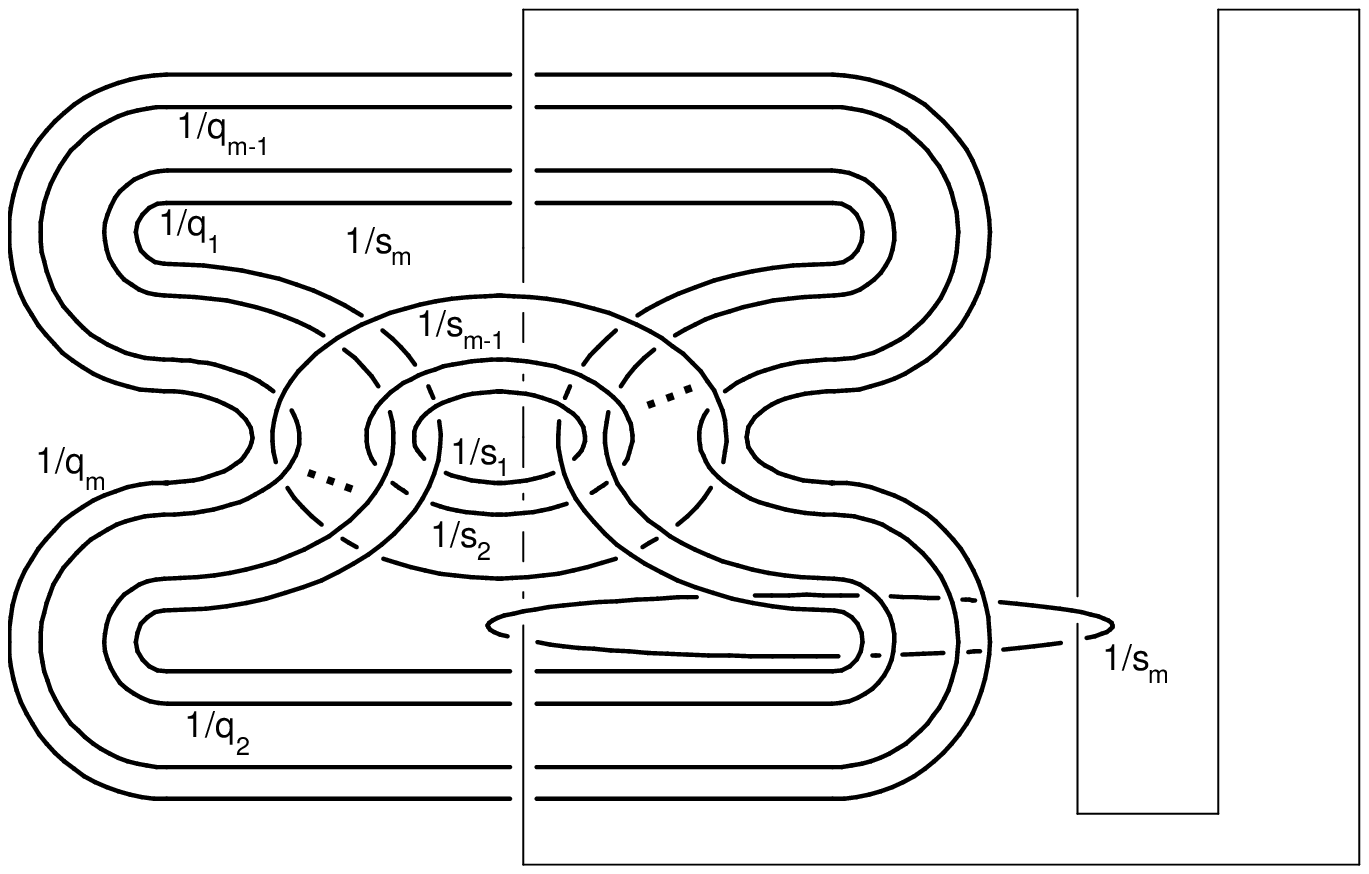}
 \end{center}
 \caption{}
 \label{Fig. 8}
\end{figure}

\begin{figure}
 \begin{center}
 \includegraphics*[totalheight=8cm]{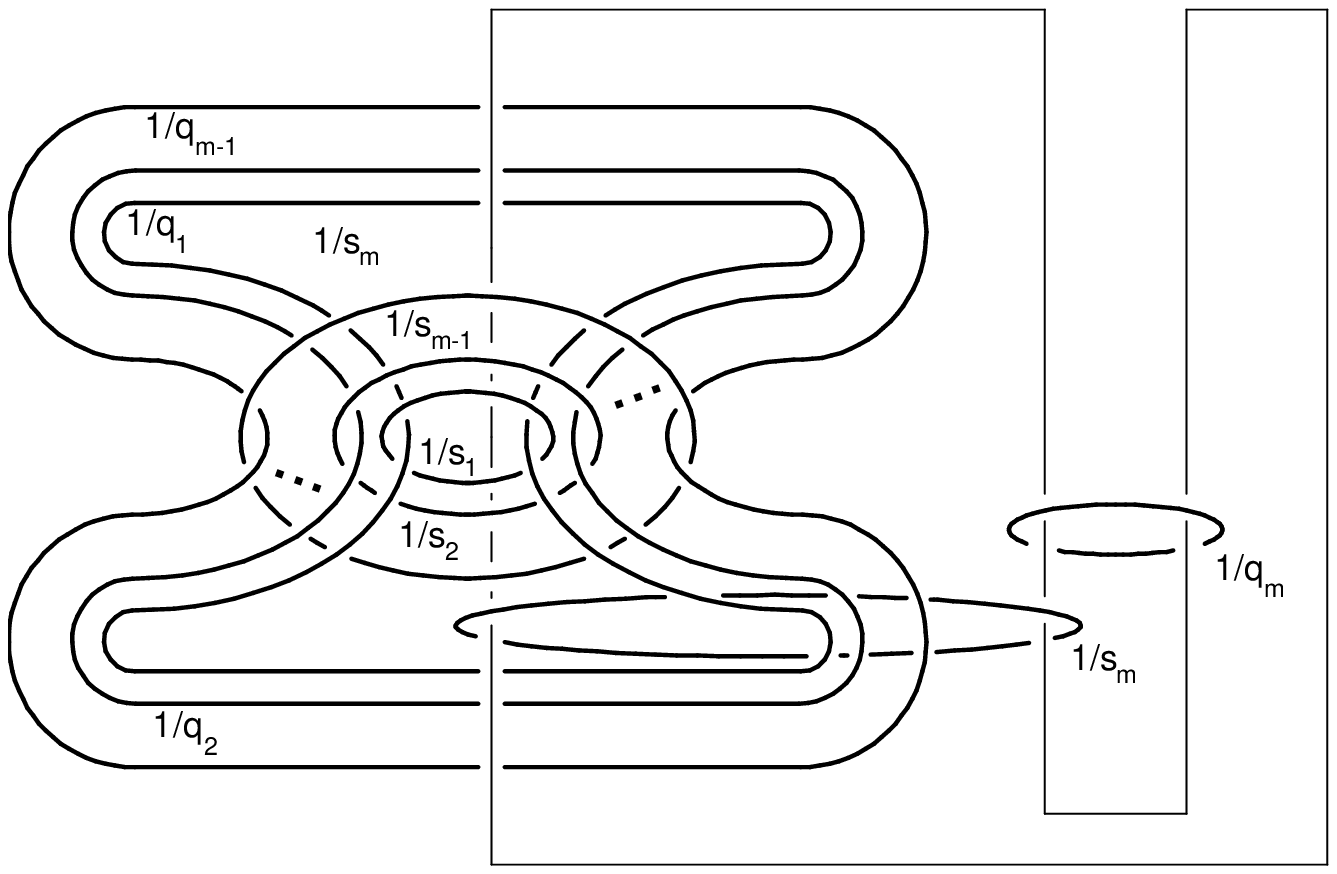}
 \end{center}
 \caption{}
 \label{Fig. 9}
\end{figure}

\begin{figure}
 \begin{center}
 \includegraphics*[totalheight=8cm]{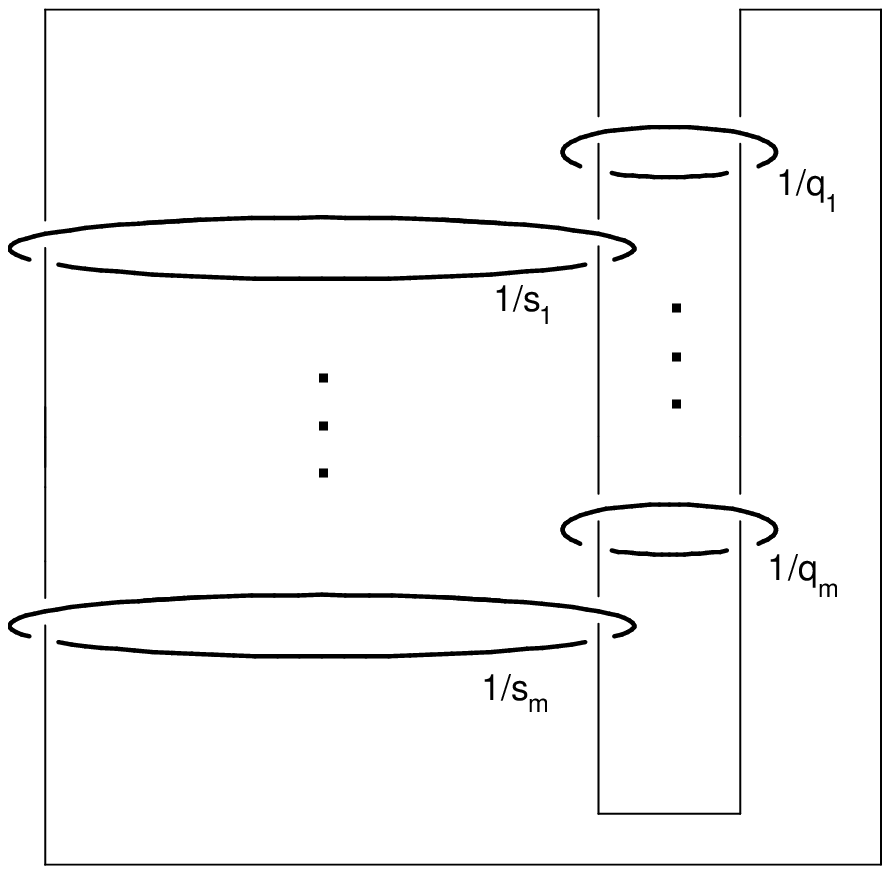}
 \end{center}
 \caption{}
 \label{Fig. 10}
\end{figure}


\begin{figure}
 \begin{center}
 \includegraphics*[totalheight=5cm]{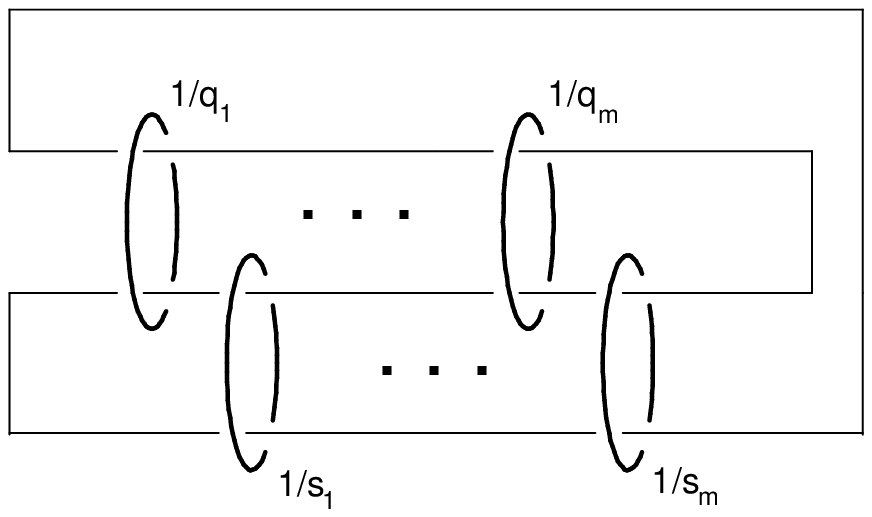}
 \end{center}
\caption{}
 \label{Fig. 11}
\end{figure}

\begin{figure}
 \begin{center}
 \includegraphics*[totalheight=5cm]{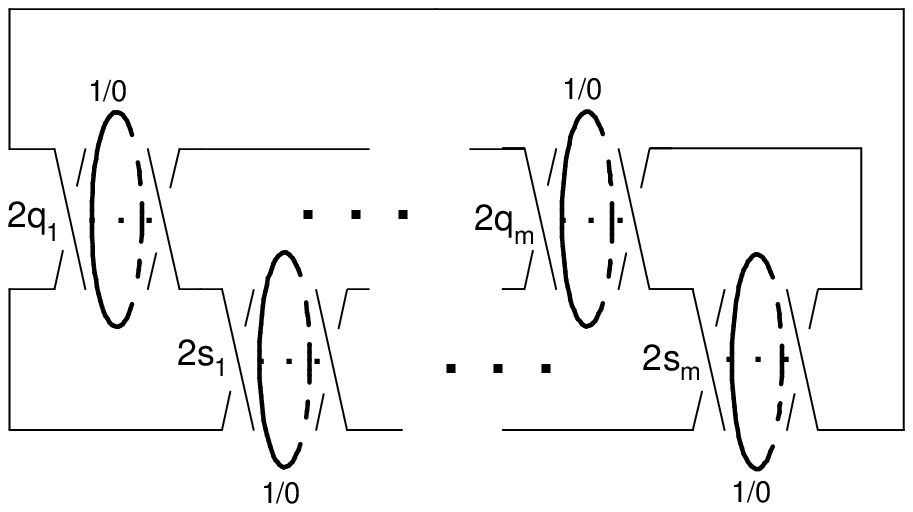}
 \end{center}
 \caption{}
 \label{Fig. 12}
\end{figure}

\begin{figure}
 \begin{center}
 \includegraphics*[totalheight=5cm]{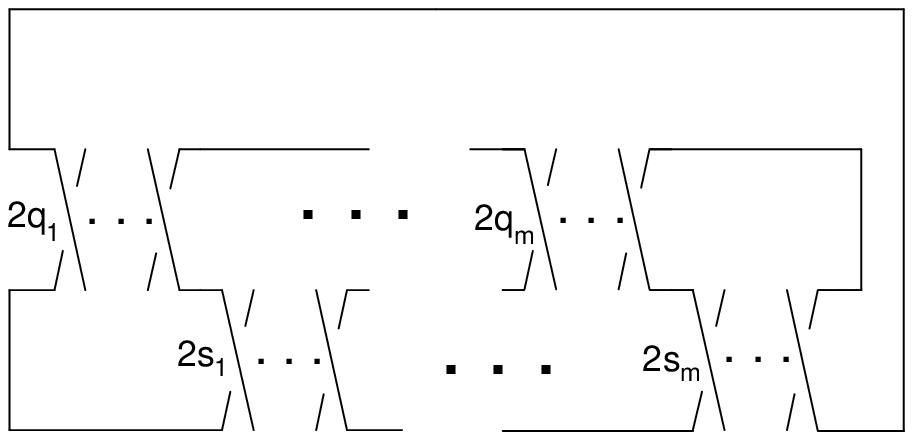}
 \end{center}
 \caption{}
 \label{Fig. 13}
\end{figure}



\vspace{15 pt} {MICHELE MULAZZANI, Department of Mathematics,
University of Bologna, I-40127 Bologna, ITALY, and C.I.R.A.M.,
Bologna, ITALY. E-mail: mulazza@dm.unibo.it}

\vspace{15 pt} {ANDREI VESNIN, The Sobolev Institute of
Mathematics, Novosibirsk, 630090 Russia. E-mail:
vesnin@math.nsk.ru}

\end{document}